\documentclass[12pt]{article}
\usepackage{amssymb,amsmath}
\usepackage{cases}
\usepackage{amsfonts}
\usepackage{color,xcolor}
\usepackage[left=2.0cm,right=2.0cm,top=2.0cm,bottom=2.0cm]{geometry}
\usepackage[colorlinks,citecolor=blue,urlcolor=blue]{hyperref}

\newtheorem{theorem}{Theorem}[section]

\newtheorem{remark}{Remark}[section]

\newcommand{\bal}{\begin{align}}
\newcommand{\bbal}{\begin{align*}}
\newcommand{\beq}{\begin{equation}}
\newcommand{\eeq}{\end{equation}}
\newcommand{\bca}{\begin{cases}}
\newcommand{\eca}{\end{cases}}
\def\div{\mathord{{\rm div}}}
\newcommand{\pa}{\partial}
\newcommand{\fr}{\frac}
\newcommand{\na}{\nabla}
\newcommand{\De}{\Delta}

\newcommand{\cd}{\cdot}
\newcommand{\ep}{\varepsilon}
\newcommand{\dd}{\mathrm{d}}

\newcommand{\R}{\mathbb{R}}

\newcommand{\les}{\lesssim}

\begin{document}
\title{A class large solution of the 3D Hall-magnetohydrodynamic equations}

\author{Jinlu Li$^{1}$\footnote{E-mail: lijinlu@gnnu.cn}, Yanghai Yu$^{2}$\footnote{E-mail: yuyanghai214@sina.com( Corresponding author)} and Weipeng Zhu$^{3}$\footnote{E-mail: mathzwp2010@163.com}\\
\small $^1$\it School of Mathematics and Computer Sciences, Gannan Normal University, Ganzhou 341000, China\\
\small $^2$\it School of Mathematics and Statistics, Anhui Normal University, Wuhu, Anhui, 241002, China\\
\small $^3$\it School of Mathematics and Information Science, Guangzhou University, Guangzhou 510006, China}

\date{\today}

\maketitle\noindent{\hrulefill}

{\bf Abstract:} In this paper, we establish the global existence to the three-dimensional incompressible
Hall-MHD equations for a class of large initial data, whose $L^{\infty}$ norms can be arbitrarily large.


\vskip0mm\noindent{\hrulefill}

\section{Introduction}\label{sec1}
This paper focuses on the following 3D incompressible Hall-magnetohydrodynamics (Hall-MHD) equations
\begin{eqnarray}\label{3D-hmhd}
        \left\{\begin{array}{ll}
          \partial_tu+u\cd\na u-\mu\De u+\na p=b\cd\na b,& x\in \R^3,t>0,\\
          \partial_tb+u\cd\na b-\nu\Delta b+\nabla\times((\na\times b)\times b)=b\cd\na u,& x\in \R^3,t>0,\\
         \div u=\div b=0,& x\in \R^3,t\geq0,\\
          (u,b)|_{t=0}=(u_0,b_0),& x\in \R^3,\end{array}\right.
        \end{eqnarray}
where $u=(u_1(t,x),u_2(t,x),u^3(t,x))\in\R^3$ and $b=(b_1(t,x),b_2(t,x),b^3(t,x))\in\R^3$ denote the divergence free velocity field and magnetic field, respectively, $p\in \R$ is the scalar pressure. $\mu$ is the viscosity and $\nu$ is the magnetic diffusivity. $\Lambda=(-\Delta)^{\fr12}$ is the Zygmund operator and the fractional power operator $\Lambda^{\gamma}$ with $0<\gamma<1$ is defined by Fourier multiplier with symbol $|\xi|^{2\gamma}$ (see e.g. \cite{Jacob 2005}), namely,
\begin{eqnarray*}
   \Lambda^{\gamma}u(x)=\mathcal{F}^{-1}|\xi|^{2\gamma}\mathcal{F}u(\xi).
\end{eqnarray*}
The parameter $\mu$ denotes the kinematic viscosity coefficient of the fluid and $\nu$ denotes the reciprocal of the magnetic Reynolds number. Comparing with the standard MHD system , the Hall-MHD system has the Hall term $\nabla\times((\na\times b)\times b)$ due to the happening of the magnetic reconnection in the case of large magnetic shear and describes many physical phenomena such as magnetic reconnection in space plasmas \cite{Forbes 1991}, star formation \cite{Balbus 2001,Wardle 2004}, neutron stars \cite{Shalybkov 1997} and also geo-dynamo \cite{Mininni 2003}. For the physical background of the magnetic reconnection and the Hall-MHD, we refer the readers to \cite{Forbes 1991,Lighthill 1960,Simakov 2008} and references therein.

Due to the physical applications and mathematical significance, there have been huge literatures on the study of the problem \eqref{3D-hmhd} by various authors. The Hall-MHD equations from two fluids model or kinetic models were derived in a mathematically rigorous way by Acheritogaray--Degond--Frouvelle--Liu \cite{Acheritogaray 2011}, where the global weak solutions in the periodic setting $\mathbb{T}^3$ by the Galerkin approximation were established. Later, Chae--Degond--Liu \cite{Chae 2014} obtained the global existence of weak solutions and the local well-posedness of classical solution to the Hall-MHD system in the whole space $\R^3$. Furthermore, they also proved that the local smooth solutions are global in time for small initial data in Sobolev space $H^s(\R^3)$ with $s>\fr52$. It is well known that the system \eqref{3D-hmhd} reduces to the classical magnetohydrodynamics (MHD) equations when the Hall term is neglected,  wherer the velocity field plays a more dominant role than the magnetic field does. Since the Hall term $\nabla\times((\na\times b)\times b)$ in the magnetic field is a quadratic term and contains the second order derivatives, its presence in the Hall-MHD system leads to the big obstacle to establish in the proof of local well-posedness. Actually, the solution of \eqref{3D-hmhd} may be ill-posedness \cite{Chae 2016} or blow-up in finite time \cite{Chae 2015} provided without this diffusion or with only fractional megnetic diffusion $\Lambda^{\alpha}b$ with $0<\alpha\leq\fr12$. The mathematical studies on \eqref{3D-hmhd} have motivated a large number of research papers concerning the low regularity local well-posedness \cite{Wu1 2018,Wan 2017}, regularity criterions \cite{He 2016,Dai 2016}, asymptotic behavior \cite{Wu2 2017,Chae 2013} and we can refer the readers to the reference therein. Since there is no global well-posedness theory for general initial data, many literatures have been devoted to the study of global existence of smooth solutions to \eqref{3D-hmhd} with small initial data.  Chae--Lee \cite{Lee 2014} gave two global-in-time existence results of the classical solutions for small initial data $(u_0,b_0)\in\dot{H}^{\fr32}(\R^3)$ or $(u_0,b_0)\in\dot{B}^{\fr12}_{2,1}(\R^3)$ who improved significantly the results in \cite{Chae 2014}. Later, Wan--Zhou \cite{Wan 2015} extended the conditions of global existence provided initial norms $||u_0||_{\dot{H}^{\fr12+\varepsilon}}+||b_0||_{\dot{H}^{\fr32}}$ are sufficiently small where $ \varepsilon\in(0,1)$. Very recently, Wan--Zhou proved that global existence of strong solution for \eqref{3D-hmhd} with the Fujita-Kato type initial data (namely, $||u_0||_{\dot{H}^{\fr12}}$ is sufficiently small), for more details, see \cite{Wan 2019}. It is also worth to mention that when $b=0$, the system \eqref{3D-hmhd} is reduced to the Navier-Stokes equations. Lei--Lin--Zhou \cite{Lei 2015} constructed a family of finite energy smooth large solutions to the Navier-Stokes equations with the initial data close to a Beltrami flow. Li--Yang--Yu \cite{Li 2019} established a class global large solution to the 2D MHD equations with damp terms whose initial energy can be arbitrarily large. Motivated by the ideas that used in \cite{Lei 2015,Li 2019}, we expect the system \eqref{3D-hmhd} with $\mu=\nu$ can generate unique global solutions for some class of large initial data.

We assume from now on that the coefficients $\mu=\nu=1$, just for simplicity. Our main result is stated as follows.
\begin{theorem}\label{the1.1} Let $U_0$ be a smooth function satisfying ${\rm{div}}U_0=0$, $\na\times U_0=\Lambda U_0$ and
\begin{eqnarray}\label{Equ1.2}
\mathrm{supp} \ \hat{U}_0(\xi)\subset\mathcal{C}\triangleq\Big\{\xi\in \R^3 \big| \ 1-\ep\leq  |\xi|\leq 1+\ep\Big\},\quad 0<\ep<\frac{2-\sqrt{2}}{2}.
\end{eqnarray}
Assume that the initial data fulfills $u_0=U_0$ and $b_0=-\na\times U_0$, then there exists a sufficiently small positive constant $\delta$, and a universal constant $C$ such that if
\begin{align}\label{condition}
C\ep^4||U_0||^2_{L^2}\Big(||\hat{U}_0||^2_{L^1}+||{U_0}||^2_{L^2}\Big)\exp\Big( C(||\hat{U}_0||_{L^1}+||\hat{U}_0||^2_{L^1})\Big)\leq \delta,
\end{align}
then the system \eqref{3D-hmhd} has a unique global solution.
\end{theorem}

\begin{remark}\label{rem1.1}
Assume that $\hat{a}_0: {\mathbb R}^3\to [0, 1]$ be a radial, non-negative, smooth function which is supported in $\mathcal{C}$ and
$\hat{a}_0\equiv 1$ for $1-\frac12\ep\leq  |\xi|\leq 1+\frac12\ep$.

Notice that
$$a_0(x)=\int_{\R^3}\cos(x\cdot\xi)\hat{a}_0(\xi)\dd \xi \quad\mbox{and}\quad \Lambda^{-1} a_0(x)=\int_{\R^3}\cos(x\cdot\xi)\frac{\hat{a}_0(\xi)}{|\xi|}\dd \xi,$$
then we have $a_0,\Lambda^{-1} a_0\in \R$ and also let $U_0=V_0+\Lambda^{-1} \na\times V_0$ with
\begin{eqnarray*}
&V_0=\ep^{-1}\Big(\log\log\frac1\ep\Big)^\frac12\na\times
\begin{pmatrix}
a_0 \\ 0 \\ 0
\end{pmatrix}
=\ep^{-1}\Big(\log\log\frac1\ep\Big)^\frac12
\begin{pmatrix}
0 \\ \pa_3a_0 \\ -\pa_2a_0
\end{pmatrix}.
\end{eqnarray*}
Here, we can show that $\mathrm{div} U_0=0$ and $\na \times U_0=\Lambda U_0$.

Moreover, we also have
\begin{eqnarray*}
\hat{U}_0=\ep^{-1}\Big(\log\log\frac1\ep\Big)^\frac12
\begin{pmatrix}
\xi^2_2+\xi^2_3 \\ -\xi_1\xi_2+i\xi_3|\xi| \\ -\xi_1\xi_3-i\xi_2|\xi|
\end{pmatrix}\frac{\hat{a}_0(\xi)}{|\xi|}.
\end{eqnarray*}
Then, direct calculations show that
\begin{align*}
||\hat{U}_0||_{L^1}\approx \Big(\log\log\frac1\ep\Big)^\frac12\quad\mbox{and}\quad||{U}_0||_{L^2}\approx \ep^{-\fr12}\Big(\log\log\frac1\ep\Big)^\frac12
\end{align*}
Thus, the left side of \eqref{condition} becomes
\begin{align*}
C\ep^2\Big(\log\log\frac1\ep\Big)^2\exp\Big(C\log\log \frac1\ep\Big).
\end{align*}
Therefore, choosing $\ep$ small enough, we deduce that the system \eqref{3D-hmhd} has a global solution.

Moreover, we also have
\begin{align*}
&||\hat{u}_0||_{L^1}\gtrsim \Big(\log\log\frac1\ep\Big)^\frac12\qquad\mbox{and} \qquad ||\hat{b}_0||_{L^1}\gtrsim \Big(\log\log\frac1\ep\Big)^\frac12.
\end{align*}
Notice that $U^1_0=-\ep^{-1}\Big(\log\log\frac1\ep\Big)^\frac12(\pa^2_{2}+\pa^2_{3})\Lambda^{-1} a_0$ and $\hat{U}^1_0\geq 0$, we can deduce that
\bbal
||U^1_0||_{L^\infty}\approx ||\hat{U}^1_0||_{L^1}\gtrsim \Big(\log\log\frac1\ep\Big)^\frac12.
\end{align*}
Since $b^1_0=-\Lambda U^1_0$, we also have $||b_0||_{L^\infty}\gtrsim \Big(\log\log\frac1\ep\Big)^\frac12$.
\end{remark}
{\bf Notations}: Let $\alpha=(\alpha_1,\alpha_2,\alpha_3)\in \mathbb{N}^3$ be a multi-index and $D^{\alpha}=\pa^{|\alpha|}/\pa^{\alpha_1}_{x_1}\pa^{\alpha_2}_{x_2}\pa^{\alpha_3}_{x_3}$ with $|\alpha|=\alpha_1+\alpha_2+\alpha_3$. For the sake of simplicity, $a\lesssim b$ means that $a\leq Cb$ for some ``harmless" positive constant $C$ which may vary from line to line. $[A,B]$ stands for the commutator operator $AB-BA$, where $A$ and $B$ are any pair of operators on some Banach space $X$. We also use the notation $||f_1,\cdots,f_n||_{X}\triangleq||f_1||_{X}+\cdots+||f_n||_{X}$.

\section{Reformulation of the System}\label{sec2}
\setcounter{equation}{0}
Denoting $\omega=\na\times u$, then we can rewrite $\eqref{3D-hmhd}_1$  as
\bal\label{3}
\pa_t\omega+u\cd\na\omega-\omega\cd\na u-\De \omega=\na\times((\na\times b)\times b).
\end{align}
It should be noted that the divergence free condition $\div b=0$ for the magnetic field will preserve all the times if $\div b_0=0$.\\
Introducing the new quantity $\Omega=\omega+b$, we can deduce from \eqref{3} and $\eqref{3D-hmhd}_2$ that
\bal\label{4}
\pa_t\Omega-\Omega\cd\na u+u\cd\na \Omega-\De \Omega=0.
\end{align}
Since $\Omega_0=0$, then we have $\Omega=0$. In other words, if $b_0=-\na\times u_0$, the relation $b=-\na\times u$ will preserve all the times.\\

Let $U$ be the solutions of the heat equations
\bal\label{1}
 \pa_tU-\De U=0, \quad  U|_{t=0}=U_0.
\end{align}
Setting $B=-\na\times U$, from \eqref{1}, we know that $(U,B)=(e^{t\Delta}U_0,e^{t\Delta}B_0)$ solve the following system
\begin{eqnarray}\label{2}
        \left\{\begin{array}{ll}
          \pa_t U-\Delta U=0,\\
          \pa_t B-\Delta B=0,\\
          \div U=\div B=0,\\
          (U,B)|_{t=0}=(U_0,B_{0}).\end{array}\right.
\end{eqnarray}

Denoting $v=u-U$ and $c=b-B$, we can reformulate the system \eqref{3D-hmhd} and \eqref{2} equivalently as
\begin{eqnarray}\label{5}
        \left\{\begin{array}{ll}
\partial_tv+v\cd\na v-c\cd \na c-\De v+\na \Big(p+\frac{|U|^2-|B|^2}{2}\Big)=f+f_1,\\
\partial_tc+v\cd\na c-c\cd\na v-\De c+\nabla\times((\na\times c)\times c)=g+g_1+g_2,\\
\div v=\div c=0,\\
(v,c)|_{t=0}=(0,0),\end{array}\right.
\end{eqnarray}
where
\bbal
&f=B\cd\na B-U\cd\na U-\na\Big(\frac{|B|^2-|U|^2}{2}\Big),
\\&f_1=B\cd\na c+c\cd\na B-U\cd\na v-v\cd\na U,
\\&g=B\cd\na U-U\cd\na B-\na\times((\na\times B)\times B),
\\&g_1=B\cd\na v+c\cd\na U-U\cd\na c-v\cd\na B,
\\&g_2=-\nabla\times((\na\times c)\times B)-\nabla\times((\na\times B)\times c).
\end{align*}

As mentioned above, although it seems difficult to control the second order derivatives in the Hall term by the diffusion term, we find that the Hall term can provide us the nice structure of the magnetic equations. Precisely speaking, both the terms $(\na\times B)\times B-(\na\times U)\times U$ and $B\cd\na U-U\cd\na B-\na\times((\na\times B)\times B)$ can generate a small quantity which is the key to construct global solutions under some assumptions of the initial data.\\
Due to $\div U=\div B=0$, we have
$$(\na \times U)\times U=U\cd\na U-\na\Big(\frac{|U|^2}{2}\Big),$$
$$(\na \times B)\times B=B\cd\na B-\na\Big(\frac{|B|^2}{2}\Big),$$
and
$$\na \times(U\times B)=B\cd\na U-U\cd\na B.$$

Notice that $B=-\na\times U=-\Lambda U$ and $\na\times B=\De U$, then we can show that
\bbal
f&=B\cd\na B-U\cd\na U-\na\Big(\frac{|B|^2-|U|^2}{2}\Big)
\\&=\De U\times B+B\times U
\\&=(\De-\mathbb{I}+2\Lambda) U\times B
\end{align*}
and
\bbal
g&=B\cd\na U-U\cd\na B-\na\times((\na\times B)\times B)
\\&=\na\times((\mathbb{I}-\De-2\Lambda)U\times B).
\end{align*}
\section{Proof of Theorem \ref{the1.1}}\label{sec3}
\setcounter{equation}{0}
In this section, we present the proof of Theorem \ref{the1.1}.\\
{\bf Proof of Theorem \ref{the1.1}}\quad   Applying $D^\alpha$ on $\eqref{5}_1$ and $\eqref{5}_2$ respectively and taking the scalar product of them with $D^\alpha v$ and $D^\alpha c$ respectively, adding them together and then summing the result over $|\alpha|\leq 3$, we get
\bal\label{z0}
\fr12\frac{\dd}{\dd t}\Big(||v||^2_{H^3}+||c||^2_{H^3}\Big)+||\na v||^2_{H^3}+||\na c||^2_{H^3}\triangleq\sum^{11}_{i=1}I_i,
\end{align}
where
\bbal
&I_1=-\sum_{0<|\alpha|\leq 3}\int_{\R^3}[D^{\alpha},v\cd] \na v\cd D^\alpha v\dd x
-\sum_{0<|\alpha|\leq 3}\int_{\R^3}[D^{\alpha},v\cd] \na c\cd D^\alpha c\dd x,
\\&I_2=\sum_{0<|\alpha|\leq 3}\int_{\R^3}[D^{\alpha},c\cd] \na c\cd D^\alpha v\dd x
+\sum_{0<|\alpha|\leq 3}\int_{\R^3}[D^{\alpha},c\cd] \na v\cd D^\alpha c\dd x,
\\&I_3=\sum_{0<|\alpha|\leq 3}\int_{\R^3}D^{\alpha}((\na\times c)\times c)\cd  D^{\alpha}(\nabla\times c)\dd x,
\\&I_4=-\sum_{0<|\alpha|\leq 3}\int_{\R^3}D^{\alpha}(U\cd \na v)\cd D^\alpha v\dd x-\sum_{0<|\alpha|\leq 3}\int_{\R^3}D^{\alpha}(U\cd \na c)\cd D^\alpha c\dd x,
\\&I_5=\sum_{0<|\alpha|\leq 3}\int_{\R^3}D^{\alpha}(B\cd \na c)\cd D^{\alpha}v\dd x+\sum_{0<|\alpha|\leq 3}\int_{\R^3}D^{\alpha}(B\cd \na v)\cd D^{\alpha}c\dd x,
\\&I_6=\sum_{0\leq|\alpha|\leq 3}\int_{\R^3}D^{\alpha}(c\cd \na B)\cd D^{\alpha}v\dd x-\sum_{0\leq|\alpha|\leq 3}\int_{\R^3}D^{\alpha}(v\cd \na B)\cd D^{\alpha}c\dd x,
\\&I_7=\sum_{0\leq|\alpha|\leq 3}\int_{\R^3}D^{\alpha}(c\cd \na U)\cd D^{\alpha}c\dd x-\sum_{0\leq|\alpha|\leq 3}\int_{\R^3}D^{\alpha}(v\cd \na U)\cd D^{\alpha}v\dd x,
\end{align*}
\bbal
&I_8=\sum_{0<|\alpha|\leq 3}\int_{\R^3}D^{\alpha}((\na\times c)\times B)\cd D^{\alpha}(\nabla\times c)\dd x,
\\&I_9=\sum_{0\leq|\alpha|\leq 3}\int_{\R^3}D^{\alpha}((\na\times B)\times c)\cd D^{\alpha}(\nabla\times c)\dd x,
\\&I_{10}=\sum_{0<|\alpha|\leq 3}\int_{\R^3}D^{\alpha}[(\De-\mathbb{I}+2\Lambda) U\times B]\cd D^{\alpha}v\dd x+\int_{\R^3}[(\De-\mathbb{I}+2\Lambda) U\times B]\cd v\dd x,
\\&I_{11}=\sum_{0\leq |\alpha|\leq 3}\int_{\R^3}D^{\alpha}[(\mathbb{I}-\De-2\Lambda)U\times B]\cd D^{\alpha}(\nabla\times c)\dd x.
\end{align*}
Next, we need to estimate the above terms one by one.

According to the commutate estimate (See \cite{Majda 2001}),
\bal\label{l0}
\sum_{|\alpha|\leq m}||[D^{\alpha},\mathbf{g}]\mathbf{f}||_{L^2}\leq C(||\mathbf{f}||_{H^{m-1}}||\na \mathbf{g}||_{L^\infty}+||\mathbf{f}||_{L^\infty}||\mathbf{g}||_{H^m}),
\end{align}
we obtain
\bal
I_1\leq&~\sum_{0<|\alpha|\leq 3}||[D^{\alpha},v\cd] \na v||_{L^2}||\na v||_{H^2}+\sum_{0<|\alpha|\leq 3}|||[D^{\alpha},v\cd]\na c||_{L^2}||\na c||_{H^2}\nonumber\\
\leq&~C||\na v||_{L^\infty}||v||_{H^3}||\na v||_{H^2}+C||v||_{H^3}||\na c||^2_{H^2}\nonumber\\
\leq&~C||v||_{H^3}\Big(||\na v||^2_{H^3}+||\na c||^2_{H^3}\Big),\label{z1}\\
I_2\leq&~\sum_{0<|\alpha|\leq 3}||[D^{\alpha},c\cd] \na c||_{L^2}||\na v||_{H^2}+\sum_{0<|\alpha|\leq 3}|||[D^{\alpha},c\cd] \na v||_{L^2}||\na c||_{H^2}\nonumber\\
\leq&~C||\na c||_{H^2}||\na v||_{H^2}||c||_{H^3}\nonumber\\
\leq&~C||c||_{H^3}\Big(||\na v||^2_{H^3}+||\na c||^2_{H^3}\Big).\label{z2}
\end{align}
Using the cancelation equality $[D^{\alpha}(\na\times c)\times c]\cd D^{\alpha}(\nabla \times c)=0$ and \eqref{l0}, we get
\bal
I_3=&~\sum_{0<|\alpha|\leq 3}\int_{\R^3}[D^{\alpha},c\times](\na\times c)\cd D^{\alpha}(\nabla\times c)\dd x\nonumber\\
\leq&~C||\na c||_{L^\infty}||c||_{H^3}||\na c||_{H^3}\nonumber\\
\leq&~ C||c||^2_{H^3}||\na c||^2_{H^3}+\frac18||\na c||^2_{H^3}.\label{z3}
\end{align}
Invoking the following calculus inequality which is just a consequence of Leibniz's formula,
\bbal
\sum_{|\alpha|\leq 3}||[D^{\alpha},\mathbf{g}]\mathbf{f}||_{L^2}\leq C(||\na \mathbf{g}||_{L^\infty}+||\na^3 \mathbf{g}||_{L^\infty})||\mathbf{f}||_{H^2},
\end{align*}
we obtain
\bal
I_4\leq&~\sum_{0<|\alpha|\leq 3}||[D^{\alpha},U\cd] \na v||_{L^2}||\na v||_{H^2}+\sum_{0<|\alpha|\leq 3}||[D^{\alpha},U\cd] \na c||_{L^2}||\na c||_{H^2}\nonumber\\
\leq&~C\Big(||\na U||_{L^\infty}+||\na^3 U||_{L^\infty}\Big)\Big(||v||^2_{H^3}+||c||^2_{H^3}\Big),\label{z4}\\
I_5\leq&~\sum_{0<|\alpha|\leq 3}||[D^{\alpha},B\cd] \na c||_{L^2}||\na v||_{H^2}+\sum_{0<|\alpha|\leq 3}||[D^{\alpha},B\cd] \na v||_{L^2}||\na c||_{H^2}\nonumber\\
\leq&~C\Big(||\na B||_{L^\infty}+||\na^3 B||_{L^\infty}\Big)\Big(||v||^2_{H^3}+||c||^2_{H^3}\Big),\label{z5}\\
I_8\leq&~\sum_{0<|\alpha|\leq 3}||[D^{\alpha},B\times] (\na\times c)||_{L^2}||\na c||_{H^3}\nonumber\\
\leq&~C\Big(||\na B||_{L^\infty}+||\na^3 B||_{L^\infty}\Big)||\na c||_{H^2}||\na c||_{H^3}\nonumber\\
\leq&~C\Big(||\na B||_{L^\infty}+||\na^3 B||_{L^\infty}\Big)^2||c||^2_{H^3}+\frac18||\na c||^2_{H^3}.\label{z6}
\end{align}
By Leibniz's formula and H\"{o}lder's inequality, one has
\bal
I_6\leq&~||c\cd \na B||_{H^3}||v||_{H^3}+||v\cd \na B||_{H^3}||c||_{H^3}\nonumber\\
\leq&~ C\Big(||\na B||_{L^\infty}+||\na^4 B||_{L^\infty}\Big)\Big(||v||^2_{H^3}+||c||^2_{H^3}\Big),\label{z7}\\
I_7\leq&~||c\cd \na U||_{H^3}||v||_{H^3}+||v\cd \na U||_{H^3}||c||_{H^3}\nonumber\\
\leq&~ C\Big(||\na U||_{L^\infty}+||\na^4 U||_{L^\infty}\Big)\Big(||v||^2_{H^3}+||c||^2_{H^3}\Big),\label{z8}\\
I_9\leq&~ ||(\na\times B)\times c||_{H^3}||\na c||_{H^3}\nonumber\\
\leq&~ C\Big(||\na B||_{L^\infty}+||\na^4 B||_{L^\infty}\Big)||c||_{H^3}||\na c||_{H^3}\nonumber\\
\leq&~C\Big(||\na B||^2_{L^\infty}+||\na^4 B||^2_{L^\infty}\Big)||c||^2_{H^3}+\frac18||\na c||^2_{H^3}.\label{z9}
\end{align}
Owing to the fact $\dot{H}^1(\R^3)\hookrightarrow L^6(\R^3)$, we deduce
\bal
I_{10}\leq&~ ||(\De-\mathbb{I}+2\Lambda) U\times B||_{H^3}||\na v||_{H^2}+||(\mathbb{I}-\De-2\Lambda)U||_{L^2}||v||_{L^6}||B||_{L^3}\nonumber\\
\leq&~C\Big(||B||_{L^\infty}+||\na^3 B||_{L^\infty}+||B||_{L^3}\Big)||(\mathbb{I}-\De-2\Lambda)U||_{H^3}||\na v||_{H^2}\nonumber\\
\leq&~C\Big(||B||^2_{L^\infty}+||\na^3 B||^2_{L^\infty}+||B||^2_{L^3}\Big)||(\mathbb{I}-\De-2\Lambda)U||^2_{H^3}+\frac18||\na v||^2_{H^3},\label{z10}\\
I_{11}\leq&~ ||(\mathbb{I}-\De-2\Lambda)U\times B||_{H^3}||\na c||_{H^3}\nonumber\\
\leq&~ C\Big(||B||_{L^\infty}+||\na^3 B||_{L^\infty}\Big)||(\mathbb{I}-\De-2\Lambda)U||_{H^3}||\na c||_{H^3}\nonumber\\
\leq&~ C\Big(||B||^2_{L^\infty}+||\na^3 B||^2_{L^\infty}\Big)||(\mathbb{I}-\De-2\Lambda)U||^2_{H^3}+\frac18||\na c||^2_{H^3}.\label{z11}
\end{align}
Putting all the estimates \eqref{z1}--\eqref{z11} together with \eqref{z0}, and using $B=-\na\times U$, we obtain
\bal\label{z12}
&\quad \frac{\dd}{\dd t}\Big(||v||^2_{H^3}+||c||^2_{H^3}\Big)+||\na v||^2_{H^3}+||\na c||^2_{H^3}\nonumber\\&\les
\Big(||v||_{H^3}+||c||_{H^3}+||c||^2_{H^3}\Big)\Big(||\na v||^2_{H^3}+||\na c||^2_{H^3}\Big)
\nonumber\\&\quad +\Big(||B,U,\na^4 B,\na^4 U||_{L^\infty}+||\na B,\na^4 B||^2_{L^\infty}\Big)\Big(||v||^2_{H^3}+||c||^2_{H^3}\Big)
\nonumber\\&\quad +\Big(||B,\na B,\na^3 B||^2_{L^\infty}+||B||^2_{L^3}\Big)||(\mathbb{I}-\De-2\Lambda)U||^2_{H^3}\nonumber\\&\les
\Big((||v||^2_{H^3}+||c||^2_{H^3})^{\fr12}+||c||^2_{H^3}\Big)\Big(||\na v||^2_{H^3}+||\na c||^2_{H^3}\Big)
\nonumber\\&\quad +\Big(\sum_{i=1}^5||\na^iU||_{L^\infty}+||\na^2 U,\na^5U||^2_{L^\infty}\Big)\Big(||v||^2_{H^3}+||c||^2_{H^3}\Big)\nonumber\\&\quad
+\Big(||\nabla U,\na^2 U,\na^4 U||^2_{L^\infty}+||\nabla U||^2_{L^3}\Big)||(\mathbb{I}-\De-2\Lambda)U||^2_{H^3}.
\end{align}
Notice that $\mathrm{supp}\ \hat{U}\subset \mathcal{C}$ and the fact $U=e^{t\Delta}U_0$, we can verify that for any integer $m\geq1$,
\bal\label{z13}
||\nabla^m U||_{L^\infty}\leq C||U||_{L^\infty}\leq Ce^{-\frac12t}||\hat{U_0}||_{L^1},\quad ||\nabla U||_{L^3}\leq C||U||_{L^2}\leq Ce^{-\frac12t}||{U_0}||_{L^2}
\end{align}
and
\bal\label{z14}
||(\mathbb{I}-\De-2\Lambda)U||_{H^m}\leq Ce^{-\frac12t}\ep^2||U_0||_{L^2},
\end{align}
where we have used the conditions \eqref{Equ1.2}.

Inserting \eqref{z13} and \eqref{z14} into \eqref{z12} yields
\bal\label{z15}
&\quad \frac{\dd}{\dd t}\Big(||v||^2_{H^3}+||c||^2_{H^3}\Big)+||\na v||^2_{H^3}+||\na c||^2_{H^3}\nonumber\\&\leq
C\Big((||v||^2_{H^3}+||c||^2_{H^3})^{\fr12}+||c||^2_{H^3}\Big)\Big(||\na v||^2_{H^3}+||\na c||^2_{H^3}\Big)
\nonumber\\&\quad+Ce^{-\fr12t}\Big(||\hat{U}_0||_{L^1}+||\hat{U}_0||^2_{L^1}\Big)\Big(||v||^2_{H^3}+||c||^2_{H^3}\Big)
\nonumber\\&\quad+Ce^{-t}\ep^4||U_0||^2_{L^2}\Big(||\hat{U_0}||^2_{L^1}+||U_0||^2_{L^2}\Big).
\end{align}
Now, we define
\bbal
\Gamma\triangleq\sup\{t\in[0,T^*):\sup_{\tau\in[0,t]}\Big(||v(\tau)||^2_{H^3}+||c(\tau)||^2_{H^3}\Big)\leq \eta\},
\end{align*}
where $\eta$ is a small enough positive constant which will be determined later on.

Assume that $\Gamma<T^*$. For all $t\in[0,\Gamma]$, we obtain from \eqref{z15} that
\bbal
\frac{\dd}{\dd t}\Big(||v||^2_{H^3}+||c||^2_{H^3}\Big) &\leq
Ce^{-\fr12t}\Big(||\hat{U}_0||_{L^1}+||\hat{U}_0||^2_{L^1}\Big)\Big(||v||^2_{H^3}+||c||^2_{H^3}\Big)
\nonumber\\&\quad+Ce^{-t}\ep^4||U_0||^2_{L^2}\Big(||\hat{U_0}||^2_{L^1}+||{U_0}||^2_{L^2}\Big),
\end{align*}
which follows from the assumption that
\bbal
||v||^2_{H^3}+||c||^2_{H^3}&\leq C\ep^4||U_0||^2_{L^2}\Big(||\hat{U}_0||^2_{L^1}+||{U_0}||^2_{L^2}\Big)\exp\Big( C(||\hat{U}_0||_{L^1}+||\hat{U}_0||^2_{L^1})\Big)\leq C\delta.
\end{align*}
Choosing $\eta=2C\delta$, thus we can get
\bbal
\sup_{\tau\in[0,t]}\Big(||v(\tau)||^2_{H^3}+||c(\tau)||^2_{H^3}\Big)&\leq \fr\eta2 \quad\mbox{for}\quad t\leq \Gamma.
\end{align*}

So if $\Gamma<T^*$, due to the continuity of the solutions, we can obtain there exists $0<\epsilon\ll1$ such that
\bbal
\sup_{\tau\in[0,t]}\Big(||v(\tau)||^2_{H^3}+||c(\tau)||^2_{H^3}\Big)&\leq \fr\eta2 \quad\mbox{for}\quad t\leq \Gamma+\epsilon<T^*,
\end{align*}
which is contradiction with the definition of $\Gamma$.

Thus, we can conclude $\Gamma=T^*$ and
\bbal
\sup_{\tau\in[0,t]}\Big(||v(\tau)||^2_{H^3}+||c(\tau)||^2_{H^3}\Big)&\leq C<\infty \quad\mbox{for all}\quad t\in(0,T^*),
\end{align*}
which implies that $T^*=+\infty$. This completes the proof of Theorem \ref{the1.1}. $\Box$

\section*{Acknowledgments} J. Li was supported by NSFC (No.11801090).

\end{document}